\documentclass[a4paper, 11pt]{article}

\usepackage[T1]{fontenc}
\usepackage[utf8]{inputenc}
\usepackage[english]{babel}
\usepackage{amsmath}
\usepackage{amssymb}
\usepackage{amsfonts}
\usepackage{amsthm}
\usepackage{enumitem}
\usepackage{moreenum}
\usepackage{mathtools}
\usepackage{csquotes}
\usepackage{hyperref}
\hypersetup{
	colorlinks = true,
	urlcolor   = blue,
	citecolor  = black,
}
\newtheorem{theorem}{Theorem}[section]

\newtheorem{proposition}[theorem]{Proposition}

\newcommand{\R}{\mathbb{R}}
\newcommand{\abs}[2][]{#1\lvert #2 #1\rvert}

\begin{document}
% body starts here --------------------
\title{On the amplitude of steady water waves with positive constant vorticity}
\date{}
\author{Evgeniy Lokharu\footnotemark[1], Erik Wahl\'en\footnotemark[1] and J\"org Weber\footnotemark[1]}
\renewcommand{\thefootnote}{\fnsymbol{footnote}}
\footnotetext[1]{Centre for Mathematical Sciences, Lund University, 221 00 Lund, Sweden}

	\maketitle

	\begin{abstract}
		For two-dimensional steady pure-gravity water waves with a unidirectional flow of constant favourable vorticity, we prove an explicit bound on the amplitude of the wave, which decays to zero as the vorticity tends to infinity. Notably, our result holds true for arbitrary water waves, that is, we do not have to restrict ourselves to periodic or solitary or symmetric waves.
	\end{abstract}
	
	\section{Introduction}
	
	In this paper we consider two-dimensional steady gravity water waves with constant vorticity, subject to the influence of only gravity, but not surface tension. Moreover, we assume that the corresponding flows are unidirectional and that the sign of the vorticity is favourable. This means that the waves we consider are \enquote{downstream} waves, that is, the underlying currents are sheared in the same direction as the that of the wave propagation. Such currents can be thought of being induced by a wind parallel to the direction of the wave propagation. In correspondence with physical intuition, these currents should flatten out the wave. Indeed, we derive an explicit bound for the amplitude of the wave profile, which also shows that the amplitudes have to become very small as the vorticity grows without limit. It is notable that the new bound is independent of the total head and the mass flux constants.
	
	We emphasize that our result holds for any water wave as described above. In particular, we do not have to assume that the wave is periodic or solitary or symmetric. Also, it does not necessarily have to be on a solution branch bifurcating from a configuration with a flat surface. Quite recently, a bound for the wave amplitude was, among other things, derived in \cite{ConsStrVarv21}, however restricted to the case of periodic and symmetric waves on such a bifurcation branch. Moreover, our proof is much simpler than the one in \cite{ConsStrVarv21} and the result also in the restricted case is sharper in the sense that it provides an $\mathcal O(\omega^{-2})$ (instead of $\mathcal O(\omega^{-1})$) decay of the amplitude as the vorticity $\omega$ tends to infinity. In the irrotational case, $\omega = 0$, some new bounds for the wave amplitude were obtained recently in \cite{Lokharu2021}.
Similar inequalities for the amplitude, but for a negative constant vorticity, are not known. In this case any bound of the form $\mathcal O(1)$ as $\omega \to -\infty$ would be of great interest.	
	
	As an important application of a priori bounds for the wave amplitude we like to mention the verification that waves of extreme forms exist. Indeed, in typical global bifurcation results (see, for example, \cite{ConsStr04,ConsStrVarv16}) one is typically left with several possible alternatives, and such an a priori bound yields that one has to approach a wave of extreme form when following the bifurcation branch. In fact, especially in the rotational case the literature proving this limiting behaviour is quite sparse and progress has been made only recently; see \cite{ConsStrVarv21,KozlovLokharu20}.
	
	Furthermore, we remark that, while we derive a bound on the wave amplitude, in our situation bounds on the slope of the water surface have already been established in \cite{StraussWheeler16}, at least for periodic or solitary waves with symmetry.

	\section{Statement of the problem} \label{s:statement}
	
	We consider the classical model for two-dimensional steady waves on water of finite depth. We neglect the effects of surface tension and consider an ideal fluid of constant (unit) density. In the corresponding moving reference frame the stationary Euler equations are given by
	\begin{subequations}\label{eqn:trav}
		\begin{align}
			\label{eqn:u}
			(u-c)u_x + vu_y & = -P_x,   \\
			\label{eqn:v}
			(u-c)v_x + vv_y & = -P_y-g, \\
			\label{eqn:incomp}
			u_x + v_y &= 0, \\
			u_y - v_x &= \omega, \label{eqn:irr}
		\end{align}
		which holds true in a two-dimensional fluid domain $D_{\eta}$, defined by the inequality
		\[
		0 < y < \eta(x).
		\]
		Here $(u,v)$ are components of the velocity field, $y = \eta(x)$ is the surface profile, $c$ is the wave speed, $P$ is the pressure, and $g$ is the gravitational constant. The corresponding boundary conditions are
		\begin{alignat}{2}
			\label{eqn:kinbot}
			v &= 0&\qquad& \text{on } y=0,\\
			\label{eqn:kintop}
			v &= (u-c)\eta_x && \text{on } y=\eta,\\
			\label{eqn:dyn}
			P &= P_{\mathrm{atm}} && \text{on } y=\eta.
		\end{alignat}
	\end{subequations}
	
	We reformulate these equations in terms of a stream function $\psi$, defined implicitly by the relations
	\[
	\psi_y = c-u, \ \ \psi_x =  v.
	\]
	This determines $\psi$ up to an additive constant, while relations \eqref{eqn:kinbot},\eqref{eqn:kintop} require $\psi$ to be constant along the boundaries. Thus, by subtracting a suitable constant we assume that
	\[
	\psi = m \text{ on } y = \eta, \quad \psi = 0 \text{ on } y = 0.
	\]
	Here $m$ is the mass flux, defined by
	\[
	m = \int_0^{\eta} (c-u) dy.
	\]
	The corresponding problem for the stream function is
	\begin{subequations}\label{stream0}
		\begin{alignat}{2}
			\label{sys:stream0:lap}
			\Delta\psi + \omega&=0 &\qquad& \text{in } D_{\eta},\\
			\label{stream0:bern}
			\tfrac 12\abs{\nabla\psi}^2 +  g y  &= Q &\quad& \text{on }y=\eta,\\
			\label{stream0:kintop} 
			\psi  &= m &\quad& \text{on }y=\eta,\\
			\label{stream0:kinbot} 
			\psi  &= 0 &\quad& \text{on }y=0.
		\end{alignat}
		In what follows we will consider the case when the flow is unidirectional, that is,%the fluid domain is free from stagnation points, that is, (SOMEWHAT STRONGER THAN FREE FROM STAGNATION POINTS??)
		\begin{equation}\label{uni}
			\psi_y > 0
		\end{equation}
		everywhere in $\overline{D_\eta}$.  \\
	\end{subequations}	
	
%	Let us denote
%	\[
%	\hat{\eta} = \sup_{x \in \R}{\eta(x)}, \ \ \check{\eta} = \inf_{x \in \R}{\eta(x)}.
%	\]
%	(DO WE REALLY NEED THIS NOTATION? SEE NOTATION CLASH BELOW.)
	Our main result now can be stated as follows.
	
	\begin{theorem}\label{thm1}
		Assume that $\omega > 0$ is a constant. Let $\psi \in C^{2}(\overline{D_\eta})$ and $\eta \in C^1(\mathbb{R})$ be an arbitrary solution to \eqref{stream0} for some constants $m$ and $Q$. Then
		\[
		\sup_{x \in \R}{\eta(x)} - \inf_{x \in \R}{\eta(x)} < 2 g \omega^{-2},
		\]
		that is, the maximum possible wave amplitude is controlled by the size of the vorticity.
	\end{theorem}
	
	Our proof is essentially based on fundamental bounds for the surface profile and the Bernoulli constant. We give now all necessary details regarding these bounds.
	
	\subsection*{Bounds on the Bernoulli constant and the surface profile}
	
	Here we state the result that was proved in \cite{Kozlov2015}. We will apply it to various formulations with different gravitational constants and so we need to state it here in a general form. Thus, we consider a problem% (PROBLEM: $\hat{\eta}$ HAS TWO DIFFERENT MEANINGS)
	\begin{subequations}\label{stream1}
		\begin{alignat}{2}
			\label{sys:stream1:lap}
			\Delta\hat{\psi} + \hat{\omega}&=0 &\qquad& \text{in } D_{\hat{\eta}},\\
			\label{stream1:bern}
			\tfrac 12\abs{\nabla\hat{\psi}}^2 +  \hat{g} \hat{y}  &= \hat{Q} &\quad& \text{on }\hat{y}=\hat{\eta},\\
			\label{stream1:kintop} 
			\hat{\psi}  &= \hat{m} &\quad& \text{on }\hat{y}=\hat{\eta},\\
			\label{stream1:kinbot} 
			\hat{\psi}  &= 0 &\quad& \text{on } \hat{y}=0,
		\end{alignat}
		under  the assumption
		\begin{equation}\label{uni1}
			\hat{\psi}_{\hat{y}} > 0
		\end{equation}
		everywhere in $\overline{D_{\hat{\eta}}}$. In \eqref{sys:stream1:lap} the vorticity $\hat{\omega}$ is a positive constant.
	\end{subequations}	
	
	This system admits  a class of stream solutions that are independent of $\hat{x}$. We can compute these solutions $\hat{\psi} = \hat{\Psi}$, %(MAYBE $\hat{\Psi}$ INSTEAD?),
	$\hat{\eta} = \hat{d}$ explicitly as 
	\[
	\hat{\Psi}(\hat{y};\hat{s}) = - \tfrac12 \hat{\omega} \hat{y}^2 + \hat{s} \hat{y}, \ \ \hat{d}(\hat{s}) = \int_0^{\hat{m}} \frac{1}{\sqrt{\hat{s}^2-2\hat{\omega} p}} dp.
	\]
	The corresponding constant $\hat{Q}$ in \eqref{stream1:bern} is then given by
	\[
	\hat{Q}(\hat{s}) = \tfrac12 \hat{s}^2 - \hat{\omega} \hat{m} + \hat{g} \hat{d}(\hat{s}).
	\]
	Note that in order to satisfy \eqref{stream1:bern}, \eqref{uni1} it is required that% (IS THIS REALLY TO SATISFY \eqref{uni1}?)
	\[
	\hat{s} > \hat{s}_0 = \sqrt{2 \hat{m} \hat{\omega}}.
	\]
	Computing the derivative with respect to $\hat{s}$, we find
	\[
	\frac{\partial \hat{Q}}{\partial \hat{s}} = \hat{s} \left(1 - \hat{g} \int_0^{\hat{m}} \frac{1}{(\hat{s}^2-2\hat{\omega} p)^{\tfrac32}} dp \right).
	\]
	Thus, there exists a unique $\hat{s} = \hat{s}_c$ for which the derivative is zero and this value can be found from the equation
	\[
	\hat{g} \int_0^{\hat{m}} \frac{1}{(\hat{s}_c^2-2\hat{\omega} p)^{\tfrac32}} dp = 1.
	\]
	Let us denote by $\hat{Q}_0=\hat{g}\sqrt{2\hat{m}/\hat{\omega}}$ and $\hat{Q}_c$ the corresponding constants $\hat{Q}(\hat{s})$ for $\hat{s} = \hat{s}_0$ and $\hat{s} = \hat{s}_c$ respectively.% (MAYBE STATE THE EXPLICIT FORMULAS?)	
	
	Assume that we are given some $\hat{q} \in (\hat{Q}_c, \hat{Q}_0)$. Then the equation
	\[
	\hat{Q}(\hat{s}) = \hat{q}
	\]
	has exactly two distinct solutions $\hat{s} = \hat{s}_-(\hat{q})$ and $\hat{s} = \hat{s}_+(\hat{q})$ with $\hat{s}_- < \hat{s}_+$. The corresponding depths are given by
	\[
	\hat{d}_+(\hat{q}) = \hat{d}(\hat{s}_-(\hat{q})), \ \ \hat{d}_-(\hat{q}) = \hat{d}(\hat{s}_+(\hat{q})).
	\]
	It is defined so that $\hat{d}_- < \hat{d}_+$.  Now we are ready to formulate a part of Theorem 1 from \cite{Kozlov2015}.
	
	\begin{proposition} \label{prop1}
		Let $\hat{\psi} \in C^2(\overline{D_{\hat{\eta}}})$ and $\hat{\eta} \in C^1(\R)$ be a non-stream solution to \eqref{stream1} for some constant $\hat{Q}$ in $\R$. Then $\hat{Q} \in (\hat{Q}_c, \hat{Q}_0)$ and the following is true:
		\begin{itemize}
			\item[(i)] $\inf \hat{\eta} \geq \hat{d}_-(\hat{Q})$;
			\item[(ii)] $\sup \hat{\eta} \geq \hat{d}_+(\hat{Q})$ and $\sup \hat{\eta} \leq \hat{d}_0:=\hat{d}(\hat{s}_0)=\sqrt{2\hat{m}/\hat{\omega}}$;
		\end{itemize}
		Here $\hat{d}_-(\hat{Q})$ and $\hat{d}_+(\hat{Q})$ are well defined because  $\hat{Q} \in (\hat{Q}_c, \hat{Q}_0)$.
	\end{proposition} 
	
	A proof of this result is given in \cite{Kozlov2015} even for non-constant vorticity functions.
	
	\section{Proof of Theorem \ref{thm1}}
	Without loss of generality we can assume that $(\psi,\eta)$ is not a stream solution. Then let us consider a scaling of variables, where the relative mass flux and the vorticity are scaled to unity. More precisely, we put
	\[
	\tilde{x} = \lambda x, \ \ \tilde{y} = \lambda y,\ \ \tilde{\eta}(\tilde{x}) = \lambda \eta(\lambda^{-1} \tilde{x}), \ \ \tilde{\psi}(\tilde{x},\tilde{y}) = m^{-1}\psi(x,y)
	\]
	where
	\begin{align*}%\label{lambda}
		\lambda = \left(\frac{\omega}{m}\right)^{\tfrac12}.
	\end{align*}
	The corresponding non-dimensional problem is	
	%\begin{subequations}\label{sys:stream}
		\begin{alignat*}{2}
			%\label{sys:stream:lap}
			\Delta\tilde{\psi} +1&=0 &\qquad& \text{in } D_{\tilde{\eta}} \coloneqq \{(\tilde{x},\tilde{y}): 0 < \tilde{y} < \tilde{\eta}\},\\
			%\label{sys:stream:bern}
			\tfrac 12\abs{\nabla\tilde{\psi}}^2 + \epsilon \tilde{y}  &= \tilde{Q} &\quad& \text{on }\tilde{y}=\tilde{\eta},\\
			%\label{sys:stream:kintop} 
			\tilde{\psi}  &= 1 &\quad& \text{on }\tilde{y}=\tilde{\eta},\\
			%\label{sys:stream:kinbot} 
			\tilde{\psi}  &= 0 &\quad& \text{on }\tilde{y}=0.
		\end{alignat*}
	%\end{subequations}	
	Here
	\[
	\epsilon = \frac{g}{m^{\tfrac12}\omega^{\tfrac32}},\quad \tilde Q=\frac{Q}{m\omega}.
	\]
	%(DEFINE $\tilde Q=Q/(m\omega)$ AS WELL?)
	We will follow the same notation system as introduced in the previous section. Thus, we define in a similar way all the quantities such as $\tilde{Q}_0, \tilde{Q}_c, \tilde{d}(\tilde{s})$ and so on. The corresponding values for the original system \eqref{stream0} we will denote with the same letters but without tildes.
	
%	Below we will assume that $\epsilon \leq 1$. Otherwise, we would have $m \leq g^2 \omega^{-3}$ and then
%	\[
%	Q_0 = g \left( \frac{2m}{\omega} \right)^{\tfrac12} \leq \sqrt{2} g^2 \omega^{-2}
%	\]
%	and by Proposition \ref{prop1} we have 
%	\[
%	\hat{\eta} \leq g^{-1} Q \leq g^{-1} Q_0 \leq \sqrt{2} g \omega^{-2}
%	\]
%	as desired. Therefore, there is no loss of generality in making this assumption.
	By Proposition \ref{prop1} we have
	\[
	\tilde{Q} < \tilde{Q}_0 = \epsilon \sqrt{2}
	\]
	and
	\[
	\tilde{d}_1\coloneqq\tilde{d}_-(\tilde{Q}_0) \leq \tilde{\eta} \leq \tilde{d}_0 = \sqrt{2}.
	\]
	%(SHOULD THIS BE $\tilde{d}_-(\tilde{Q}_0)$?)
	Our aim is to show that
	\begin{equation}\label{amp}
		\tilde{d}_0 - \tilde{d}_1 < 2 \epsilon.
	\end{equation}
	If this is true, then
	\[
	\sup\eta-\inf\eta = \lambda^{-1} (\sup \tilde{\eta} - \inf \tilde{\eta}) \leq \lambda^{-1} (\tilde{d}_0 - \tilde{d}_1) < 2 \lambda^{-1} \epsilon = 2 g \omega^{-2}
	\]
	as desired. Here we used the inequality $\tilde{d}_-(\tilde{Q}) > \tilde{d}_1$, which is valid since $\tilde{Q} \in (\tilde{Q}_c,\tilde{Q}_0)$.

	 In order to verify \eqref{amp}, we consider two cases: First, assume that $\epsilon\ge\sqrt2/2$. Then, clearly, 
	 \[\tilde d_0-\tilde d_1=\sqrt2-\tilde d_1<\sqrt2\le2\epsilon\]
	 and \eqref{amp} holds.
	 
	 If on the other hand $\epsilon<\sqrt2/2$, we first note that
	 \[(\sqrt2,\infty)\ni\tilde s\mapsto \tilde d(\tilde s)=\int_0^1 \frac{1}{\sqrt{\tilde s^2-2p} } \, dp=\tilde s-\sqrt{\tilde s^2-2}\in(0,\sqrt2)\]
	 is bijective with inverse
	 \[\tilde s(\tilde d)=\frac{\tilde d}{2}+\frac{1}{\tilde d}, \quad \tilde d\in(0,\sqrt{2}).\]
	 Let us correspondingly express $\tilde Q$ as a function of $\tilde d$:
	 \[
	 \tilde Q(\tilde d)=\frac12 \tilde s(\tilde d)^2-1+\epsilon \tilde d=\frac{\tilde d^2}{8}+\frac{1}{2\tilde d^2}-\frac12+\epsilon \tilde d,\quad\tilde d\in(0,\sqrt{2}). 
	 \]
%	 where
%	 \[
%	 \tilde d(\tilde s)=\int_0^1 \frac{1}{\sqrt{\tilde s^2-2p} } \, dp=\tilde s-\sqrt{\tilde s^2-2}, \quad \tilde s \ge \sqrt{2}. 
%	 \]
%	 We can in fact invert the last relation by moving terms and squaring:
%	 \[
%	\tilde s^2-2=(\tilde s-\tilde d)^2=\tilde s^2-2\tilde s \tilde d+\tilde d^2
%	 \]
%	and hence
%	\[
%	\tilde s=\frac{\tilde d}{2}+\frac{1}{\tilde d}, \quad \tilde d \le \sqrt{2}.
%	\]
%	Inserting this into the formula for $\tilde Q$, we get
%	\[
%	\tilde Q(\tilde d)=\frac{\tilde d^2}{8}+\frac{1}{2\tilde d^2}-\frac12+\epsilon \tilde d.
%	\]
	Therefore,
	\[
		\tilde Q(\sqrt{2}-\delta)-\epsilon\sqrt{2}	=-\epsilon \delta+\frac12 \delta^2+\frac{\delta^3(4\sqrt{2}-3\delta)}{8(\sqrt{2}-\delta)^2}>-\epsilon \delta+\frac12 \delta^2,\quad\delta\in(0,\sqrt2).
	\]
	Since $(0,2\epsilon)\subset(0,\sqrt2)$, it follows that
	$\tilde Q(\sqrt{2}-\delta)-\epsilon\sqrt{2}$ has a root $\delta\in(0,2\epsilon)$. This root $\delta$ has to satisfy $\sqrt2-\delta=\tilde d_1$, that is, $\tilde d_0-\tilde d_1=\delta$. Thus, \eqref{amp} also holds in this case and the proof of Theorem \ref{thm1} is complete.

%	
%	We're interested in finding $\tilde d_1<\tilde d_0\coloneqq \sqrt{2}$ such that $\tilde Q(\tilde d_1)=\sqrt{2}\epsilon$.
%	If we insert $\tilde d_1=\sqrt{2}-\delta$ in $\tilde Q(\tilde d)-\sqrt{2}\epsilon$ and expand, we get
%	\[
%	\tilde Q(\sqrt{2}-\delta)-\sqrt{2}\epsilon=\frac14-\frac{\sqrt{2}}{4}\delta+\frac{\delta^2}{8}+\frac14+\frac{\sqrt{2}}{4}\delta+\frac{3}{8}\delta^2-\frac12-\epsilon \delta
%	+\frac{\delta^3(4\sqrt{2}-3\delta)}{8(\sqrt{2}-\delta)^2}
%	=-\epsilon \delta+\frac12 \delta^2 +O(\delta^3),
%	\]
%	where  $O(\delta^3)>0$ if $\delta<4\sqrt{2}/3\approx 1.89$. Thus, there is a root $0<\delta < 2\epsilon$, if $\epsilon < 2\sqrt{2}/3$. This estimate is bad when $\epsilon>1/\sqrt{2}$, since then we get $2\epsilon>\sqrt{2}$ and clearly $\delta <\sqrt{2}$. However, when $\epsilon\to 0$ it's optimal at order $\epsilon$. Since $\sqrt{2}<4\sqrt{2}/3$, we can simply take 
%	$\delta<\min\{2\epsilon, \sqrt{2}\}$. Alternatively, we can leave it as $2\epsilon$ for simplicity. This is also correct, since $\min\{2\epsilon, \sqrt{2}\}\le 2\epsilon$.
	
	\paragraph{Declaration of interests} The authors report no conflict of interest.
	
	\paragraph{Acknowledgement} We acknowledge the support of the Swedish Research Council (grant no 2020-00440).
	
	\bibliographystyle{amsplain}
	\bibliography{bibliography}
\end{document}